\documentclass[11pt,a4wide]{article}

\usepackage{amsfonts, amsmath, amssymb, amsxtra, stmaryrd}
\usepackage{amsthm}

\usepackage[english]{babel}

\usepackage{enumerate}
\usepackage{comment}
\usepackage{mathrsfs}
\usepackage{graphicx}
\usepackage{caption,subcaption, float}
\usepackage{mathtools}
\usepackage{lscape}
\usepackage{rotating}
\usepackage{array}
\usepackage{tikz}

\usepackage{hyperref}

\numberwithin{subsection}{section}
\numberwithin{equation}{section}
\theoremstyle{plain}
\newtheorem{satz}{Theorem}[section]
\newtheorem{thm}[satz]{Theorem}

\newtheorem*{theorem*}{Theorem}
\newtheorem{lemma}[satz]{Lemma}
\newtheorem{prop}[satz]{Proposition}

\newtheorem{cor}[satz]{Corollary}

\theoremstyle{definition}

\newtheorem{remark}[satz]{Remark}

\allowdisplaybreaks

\newcommand{\ann}{\mathrm{Ann}}

\newcommand{\g}{\mathfrak{g}}
\newcommand{\h}{\mathfrak{h}}

\newcommand{\MF}{\mathbb{F}}

\newcommand{\E}{\mathcal{E}}
\newcommand{\I}{\mathcal{I}}

\newcommand{\F}{\mathbb{F}}
\newcommand{\K}{\mathbb{K}}

\newcommand{\rad}{\mathop{\mathrm{rad}}}

\begin{document}

\title{A geometric characterization of the finitary special linear and unitary Lie algebras }
\author{Hans Cuypers\footnote{corresponding author}, Marc Oostendorp\\
Department of Mathematics and Computer Science\\
Eindhoven University of Technology\\
P.O. Box 513, 5600 MB Eindhoven\\
The Netherlands\\
email: f.g.m.t.cuypers@tue.nl}

\maketitle

\begin{abstract}
An extremal element $x$ in a Lie algebra $\g$ is an element for which
the space $[x, [x, g]]$ is contained in the linear span of $x$. Long root
elements in classical Lie algebras are examples of extremal elements.
Lie algebras generated by extremal elements lead to geometries with as points
the $1$-spaces generate by extremal elements an as lines the $2$-spaces
whose non-zero elements are pairwise commuting extremal elements. 

In this paper we show that the finitary special linear Lie algebras can be characterized by their extremal geometry.
Moreover, we also show that the finitary special unitary Lie algebras
can be characterized by the fact that their geometry has no lines, but
that after extending the field quadratically, the geometry is that of a special
linear Lie algebra.
\end{abstract}

\newpage

\section{Introduction}
  An \emph{extremal} element $x$ in a Lie algebra $\mathfrak{g}$ is an element for which the space $[x,[x,\mathfrak{g}]]$ is contained in the linear span of $x$ (plus some extra conditions in case the characteristic of $\F$ is $2$). Such elemnts are called \emph{pure}, if $[x,[x,\mathfrak{g}]]$ is $1$-dimensional.
  Long root elements in classical Lie algebras are examples of extremal elements.

  If  $E$ is the set of   extremal elements of a Lie algebra $\g$,
  then we can define a geometry on the set of extremal points $\E$, which is the set of $1$-dimensional subspaces of $\g$ spanned by extremal elements.
  The lines of this geometry are those $2$-dimensional subspaces (identified with  the collection of extremal points in it) of $\g$ all whose elements
  are pairwise commuting extremal elements.
  This geometry is called the {\em extremal geometry} of $\g$ and is denoted by $\Gamma(\g)$.

  In \cite{CSUW01} Cohen et al. started the study of simple Lie algebras generated by their extremal elements, in order to  obtain
  a characterization of the classical Lie algebras.

  The main result of Cohen and Ivanyos  \cite{CI06,CI07} states that for a finite dimensional simple Lie algebra $\g$ generated by its pure extremal elements
  the extremal geometry $\Gamma(\g)$ either contains no lines or is a so-called root shadow space of a spherical building.
  Subsequently, in combined work of Fleischmann, Roberts, Shpectorov and the first author of this paper \cite{CRS,CF18},
  it has been shown that in case  $\Gamma(\g)$ is a root shadow space of a spherical building of rank at least $3$, the isomorphism type
  of $\g$ is uniquely determined by its extremal geometry.
  In particular, classical Lie algebras with an extremal geometry
  of finite rank at least $3$ are characterized by this extremal geometry.

  In this paper we continue the study of the relation between a simple Lie algebra generated by its pure extremal elements and its extremal geometry.
  In particular, we focus on (variations of) special linear and unitary  Lie algebras, both finite and infinite dimensional.

  Let $V$ be a 
  vector space over a field $\F$, then for each  $0\neq v\in V$ and $0\neq \phi\in V^*$, the dual of $V$, with $\phi(v)=0$, the linear map
  $$t_{v,\phi}:V\rightarrow V$$
  defined by $t_{v,\phi}(w)=\phi(w)v$ for all $w\in V$
  is an {\em infinitesimal transvection} in the general linear
  Lie algebra
  $\mathfrak{gl}(V)$. It is a pure extremal element, see Section \ref{section:extremal}.

  The infinitesimal transvections are finitary
  linear maps on $V$ (i.e., their kernel
  has finite codimension in $V$) and traceless,
  hence they are contained in, and actually generate,  $\mathfrak{fsl}(V)$, the  Lie algebra of all traceless
  finitary linear maps on $V$.

  Now suppose  $\Pi$ to be a subspace of $V^*$, then we can consider
  the Lie subalgebra
  $$\mathfrak{fsl}(V,\Pi):=\langle t_{v,\phi}\mid
  0\neq v\in V,\ 0\neq \phi\in \Pi\ \mathrm{and}\ \phi(v)=0\rangle$$ 
  of $\mathfrak{fsl}(V)$.

  The Lie algebras $\mathfrak{fsl}(V,\Pi)$ are
  simple up to a center, provided that
  the annihilator in $V$ of $\Pi$, i.e $\ann_V(\Pi)=\{v\in V\mid \phi(v)=0\ \mathrm{for\ all}\ \phi\in \Pi\}$, is trivial. 
  
  The extremal geometry of the Lie algebra
  $\mathfrak{fsl}(V,\Pi)$ is isomorphic
  to the geometry $\Gamma(V,\Pi)$ whose points are the incident point-hyperplane pairs
  $(p,H)$ of $\mathbb{P}(V)$ with $H$ the kernel of
  an element from $\Pi$.
  A line of $\Gamma(V,\Pi)$ consists
  of point-hyperplane pairs $(p,H)$ where $p$ is running
  over a line of $\mathbb{P}(V)$ contained in $H$,
  or, dually, where $H$ is running over all hyperplanes
  containing a codimension $2$ space on $p$.

  Our first result characterizes
  the Lie algebras $\mathfrak{fsl}(V,\Pi)$ by their extremal geometry:

  \begin{thm}\label{slthm}
    Let $V$ be a vector space over a field $\F$ of dimension at least $3$
    and $\Pi$ a subspace of $V^*$ with $\ann_V(\Pi)=\{0\}$.
    
    Suppose $\g$ is a simple Lie algebra generated by its set of pure extremal elements  with extremal geometry
    $\Gamma(\g)$ isomorphic to $\Gamma(V,\Pi)$.

    Then $\g$ is isomorphic to
    $\mathfrak{fsl}(V,\Pi)$ modulo its center.
    \end{thm}

  We will also consider several other subalgebras of $\mathfrak{fsl}(V,V^*)$ generated by infinitesimal transvections.
  If $V$ is equipped with a non-degenerate symplectic form $f$, then the infinitesimal transvections $t_{v,\phi}$
  with $0\neq v\in V $ and $\phi\in V^*$ given by $\phi(w)=f(v,w)$
  are contained in and generate the finitary symplectic Lie algebra $\mathfrak{fsp}(V,f)$, see \cite{CFsymp17}.
  The extremal geometry of this Lie algebra consist then of all point-hyperplane pairs $(p,H)$,
  where $H$ is the hyperplane of all points perpendicular to $p$ with respect to $f$.
  Obviously, this geometry does not contain lines.

  The main result of \cite{CFsymp17} can be stated as follows.

  \begin{thm}[Cuypers, Fleischmann]\label{sympthm}
    Let $\g$ be a simple Lie algebra over a field $\F$ of characteristic $\neq 2$ generated by its pure extremal elements.
    If the extremal geometry of $\g$ does not contain lines, then either
    $\g$ is isomorphic to $\mathfrak{fsp}(V,f)$ for some non-degenerate symplectic space  $(V,f)$, or
    there exists a degree $2$ extension $\K$ of $\F$ such that the extremal geometry of $\g\otimes_\F\K$ does contain lines.
  \end{thm}

  The latter occurs, for example, in the following situation.
  If $V$ is a $\K$-vector space  equipped with a non-degenerate skew-Hermitian form $h$, with respect to some field automorphism $\sigma$
  of order $2$, then the infinitesimal transvections $t_{v,\phi}$
  with $0\neq v\in V $ and $\phi\in V^*$ given by $\phi(w)=h(v,w)$ (assuming $h$ to be linear in the first coordinate)
  are contained in and generate the finitary special unitary Lie algebra $\mathfrak{fsu}(V,h)$ defined over the subfield $\F$ of $\K$ fixed by
  $\sigma$, provided $h$ is trace-valued (see Section \ref{section:transvections}).
  The extremal geometry of this Lie algebra does not contain lines.
  However, in this case the subalgebra of $\mathfrak{fsl}(V)$ (so, allowing  scalars from $\K$) generated by these transvections
  will be the algebra $\mathfrak{fsl}(V,\Pi)$, where $\Pi$ is the subspace generated by the elements $\phi\in V^*$ for which there is a vector $v\in V$  with $\phi(v)=0$ and $\phi(w)=h(w,v)$ for all $w\in V$.
In particular, its extremal geometry will be isomorphic with $\Gamma(V,\Pi)$.
  This latter situation is characterized by the following result.

  \begin{thm}\label{uthm}
Suppose $\g$ is a simple Lie algebra over a field $\F$ generated by pure extremal elements 
whose extremal geometry does not contain lines.

Assume there exists a quadratic Galois extension $\K$ of $\F$ such that the extremal geometry of $\g\otimes_\F\K$
is isomorphic to $\Gamma(V,\Pi)$ for some $\K$-vector space $V$ and subspace $\Pi$ of $V^*$
with $\ann_V(\Pi)=\{0\}$.

Then there is a trace-valued skew-Hermitian form $h$ on $V$ such that $\g$ is isomorphic to $\mathfrak{pfsu}(V,h)$. 
  \end{thm}

  The above theorem characterizes a class of simple Lie algebras generated by extremal elements for which the extremal geometry contains no lines. In recent work by
  Jeroen Meulewaeter and the first author, see \cite{CM19}, it has been shown that in characteristic different from $2,3$ such finite dimensional  Lie algebras are either symplectic, as in Theorem \ref{sympthm}, or the extremal geometry carries the structure of a Moufang polar space or Moufang set, which can be obtained as the set of fixed points of an involution acting on a root shadow space of a spherical building of rank at least $2$.
  Theorem \ref{uthm} covers the case where the root shadow space
 is obtain from a building of type $\mathrm{A}_n$.
  
  This paper is organized as follows. In Section \ref{section:extremal}
  we provide information and results on extremal elements and their geometry. In Section \ref{section:transvections} we discuss the various examples of Lie algebras appearing in the above results. Then Section \ref{section:sl} is devoted to a proof of Theorem \ref{slthm}, while the final section, Section \ref{section:u}, provides a proof of Theorem \ref{uthm}.
  
\section{Extremal elements in Lie algebras}
\label{section:extremal}

Let $\mathfrak{g}$ be a Lie algebra over the field $\mathbb{F}$  and with
Lie bracket $[\cdot,\cdot]$. 
An {\em extremal} element of $\mathfrak{g}$ 
is a nonzero element $x\in \mathfrak{g}$
with the property that there exists a map $g_x:\mathfrak{g}\rightarrow \mathbb{F}$, the \emph{extremal form} at $x$, such that
for all $y\in\mathfrak{g}$ we have
\begin{equation}[x,[x,y]]=2g_x(y)x.\label{extr}
\end{equation}
Moreover, for all $y,z\in \mathfrak{g}$ we have
\begin{equation} \begin{aligned} \label{P1}
\big[[x,y],[x,z]\big]=g_x\big([y,z]\big)x+g_x(z)[x,y]-g_x(y)[x,z]
 \end{aligned} \end{equation}
and
\begin{equation} \begin{aligned} \label{P2}
\big[x,[y,[x,z]]\big]=g_x\big([y,z]\big)x-g_x(z)[x,y]-g_x(y)[x,z]
 \end{aligned} \end{equation}
for every $y,z\in \mathfrak{g}$.\\
The last two identities are called the {\em Premet identities}.
If the characteristic of $\mathbb{F}$ is not $2$, then the Premet identities follow from Equation \ref{extr}. See \cite{CI06}.

As a consequence, $\big[ x,[x,\mathfrak{g}]\big]\subseteq \mathbb{F} x$ for an extremal $x\in \mathfrak{g}$. 
We call  $x\in \mathfrak{g}$ a {\em sandwich} if $[x,[x,y]]=0$ and $[x,[y,[x,z]]]=0$ for every $y,z\in \mathfrak{g}$. So, a sandwich is an element $x$ for which the extremal form $g_x$ can be chosen to be identically zero. 
We introduce the convention that {\em $g_x$ is identically zero whenever $x$ is a sandwich in $\mathfrak{g}$}.
An extremal element is called {\em pure} if it is not a sandwich.

We denote the set of extremal elements of a Lie algebra by $E(\mathfrak{g})$ or, if $\mathfrak{g}$   is clear from the context, by $E$.
Accordingly, we denote the set  $\{ \mathbb{F} x| x\in E(\mathfrak{g})\}$ of {\em extremal points} in the projective space on $\mathfrak{g}$ by $\mathcal{E}(\mathfrak{g})$ or $\mathcal{E}$.

We assume that $\mathfrak{g}$ is  generated by its set of extremal elements.

We recall some properties of extremal elements:

\begin{prop}\cite[Proposition 20]{CI06}\label{linear span}
The Lie algebra $\g$ is linearly spanned by its extremal elements.
\end{prop}

\begin{prop}\cite[Proposition 20]{CI06}\label{extremalform}
There is an associative symmetric bilinear form $g:\mathfrak{g}\times \mathfrak{g}\rightarrow\mathbb{F}$,
such that for all $x\in E$ and $y\in \mathfrak{g}$ we have
$$[x,[x,y]]=2g(x,y)x.$$
\end{prop}

The form $g$ is called the {\em extremal form} on $\mathfrak{g}$.
As the form $g$ is associative, its radical $\rad(g)$ is an ideal in $\mathfrak{g}$. 
Notice that by our choice to set $g_x=0$ for sandwiches, all sandwiches are in the radical of $g$.

\medskip

\begin{prop}\label{simple iff pure}
  Suppose $\g$ is a simple Lie algebra over the field $\mathbb{F}$ generated by its extremal elements.

  If $\g$ contains a pure extremal element, then all extremal elements are pure.
\end{prop}

\begin{proof}
  If $\g$ is simple then either $\rad(g)=\{0\}$ and there are no sandwiches,
  of $\rad(g)=\g$ and all extremal elements are in the radical of $g$.
  In the latter case all extremal elements are sandwiches.
\end{proof}


\begin{prop}\cite{CSUW01,CI06}\label{auto}
Let $x$ be  a pure extremal element of $E$. Then for each $\lambda\in \MF$ the map
$$\mathrm{exp}(x,\lambda):\g\rightarrow \g,$$
defined by
$$\mathrm{exp}(x,\lambda)y=y+\lambda[x,y]+\lambda^2g(x,y)x$$
for all $y\in \g$, is an automorphism of $\g$.
\end{prop}

If $x\in E$, then we denote by $\mathrm{Exp}(\langle x\rangle)$ the subgroup $\{\exp(x,\lambda)\mid\lambda\in \F\}$ of $\mathrm{Aut}(\g)$.
Notice that this definition is independent from the choice of $x$ in $\langle x\rangle$ as $\exp(\lambda x,1)=\exp(x,\lambda)$ for $0\neq \lambda\in \F$.

\medskip

\begin{prop}\cite{CSUW01,CI06}\label{2generators}
For pure $x,y\in E$  we have one of the following:

\begin{enumerate}[\rm (a)]
\item $\mathbb{F}x=\mathbb{F}y$;
\item $[x,y]=0$ and $\lambda x+\mu y\in E\cup \{0\}$ for all $\lambda,\mu\in \mathbb{F}$;
\item $[x,y]=0$ and $\lambda x+\mu y\in E$ only if $\lambda=0$ or $\mu=0$;
\item $z:=[x,y]\in E$, and $x,z$ and $y,z$ are as in case {\rm (b)};
\item $g(x,y)\neq 0$ and $\langle x,y\rangle $ is isomorphic to $\mathfrak{sl}_2(\mathbb{F})$.
\end{enumerate}

\end{prop}

An {\em extremal line} in $\mathfrak{g}$ is a 2-dimensional subspace of $\mathfrak{g}$
such that all its elements are extremal and pairwise commuting.
We identify an extremal line also with the set of extremal points contained in it.
Two linearly independent elements on an extremal line are as in case (b) of the above Proposition \ref{2generators}.

The {\em extremal geometry} $\Gamma(\g)$ of $\g$ 
is the point-line geometry with as
point set $\E$ and as lines
the extremal lines.

The $\mathfrak{sl}_2$ graph $\Gamma_{\mathfrak{sl}_2}$ of $\g$
is the graph with as vertices the extremal points and two vertices
$\langle x\rangle,\langle y\rangle$, where $x,y\in E$, adjacent
if and only if $g(x,y)\neq 0$.

The main result of Cohen and Ivanyos \cite{CI06,CI07}
is the following.

\begin{thm}[Cohen and Ivanyos]\label{ci}

  Let $\g$ be a Lie algebra generated by its pure extremal elements.
  If the extremal form $g$ is non-degenerate, then $\g$ is a direct product of Lie algebras $\mathfrak{h}$ generated by
  connected components of $\Gamma_{\mathfrak{sl}_2}$; the extremal geometry $\Gamma(\h)$ is
  either a root shadow space of a spherical building
  or does not contain lines.
\end{thm}

This result together with the following theorem,
which combines the results of \cite{CF18} and \cite{CRS}, provides us with a characterization of
most of the classical Lie algebras.

\begin{thm}[Cuypers, Fleischmann, Roberts, Shpectorov] 
   Let $\g$ be a 
  simple Lie algebra generated by its pure extremal elements.
  If the extremal geometry $\Gamma(\g)$ is
  a root shadow space of a spherical building
   of rank at least $3$, or $\Gamma(\g)$ is of type $\mathrm{A}_{2,\{1,2\}}$, then the isomorphism
   type of $\g$ is determined by its extremal geometry
   $\Gamma(\g)$.
\end{thm}

The cases not covered by the above theorem are root shadow spaces of type $\mathrm{G}_{2,2}$ and the cases where the extremal geometry has no lines.
We will now be concerned with the latter by studying Lie algebras over extension fields.


Suppose that $\g$ is a simple Lie algebra over a field $\F$
generated by its set
$E$ of pure extremal elements.

If $\K$ is a field extension of $\F$,
then the Lie product on  $\hat{\g}:=\g\otimes_{\F}\K$
is defined by

$$[x\otimes \lambda,y\otimes\mu]=[x,y]\otimes\lambda\mu,$$
for all $x,y\in \g$ and $\lambda,\mu \in \K$.

The elements $x\otimes \lambda$ with $x\in E$ and $\lambda\in \K$ are easily checked to be extremal.
Moreover, they generate $\hat{\g}$.

The extremal form $\hat{g}$ on $\hat{\g}$ satisfies
$$\hat{g}(x\otimes\lambda,y\otimes \mu)=g(x,y)\lambda\mu$$
for all $x,y\in E$ and $\lambda,\mu\in \K$.

Now we can formulate the main result from \cite{CFsymp17}:

\begin{thm}[Cuypers and Fleischmann]\label{thmsymp}
  Let $\g$ be a simple Lie algebra over a field $\F$ of characteristic not $2$
  generated by its set
  $E$ of pure extremal elements.
  Suppose that the extremal geometry $\Gamma(\g)$ does not contain lines.
  Then either $\g$ is isomorphic
  to a finitary symplectic Lie algebra $\mathfrak{fsp}(V,f)$ for some non-degenerate symplectic space $(V,f)$,
  or there is a degree $2$ extension $\K$ of $\F$
  such that $\g\otimes_{\F}\K$ is a Lie algebra
  generated by pure extremal elements and containing
  extremal lines.
\end{thm}

The above results imply that, at least if the characteristic is not $2$,
the study
of  Lie algebras generated by pure extremal
elements can be restricted
to those for which the extremal geometry does contain lines and their subalgebras over a subfield of index $2$.

In the next section we discuss how this situation arises for special linear Lie algebras and their unitary subalgebras. 
But first we will show that we may assume the Lie algebra $\hat{\mathfrak{g}}$ to be simple.

So, let $\g$ be a simple Lie algebra over a field $\F$ of characteristic not $2$.
Assume $\g$ is generated by its set $E$ of pure extremal elements.
Then the extremal form $g$ is non-degenerate and $\g$ contains no sandwiches, see Proposition \ref{simple iff pure}

Assume $\K$ to be a Galois extension of $\F$ of degree $2$ and 
let $\sigma$ be the field automorphism 
of $\K$ of order $2$ fixing precisely
the elements of $\F$.
Then $\sigma$ induces an automorphism, also denoted by  ${\sigma}$
of $\hat{\g}$ by mapping each
$x\otimes \lambda$, where $x\in \g$ and $\lambda\in \K$ to $(x\otimes\lambda)^{{\sigma}}:=x\otimes \lambda^{\sigma}$.

The elements fixed by ${\sigma}$ are precisely all
element $x\otimes \lambda$, where $x\in\g$ and $\lambda\in \F$,
forming a Lie subalgebra over $\F$ isomorphic and identified with $\g$.

The extremal form on $\hat{\g}$ is denoted by $\hat{g}$.

\begin{lemma}\label{radical trivial}
  The radical of the extremal form $\hat{g}$ is trivial.
\end{lemma}

\begin{proof}
  Suppose $r\neq 0$ is an element in the radical of $\hat{g}$
  and $r=x_1\otimes \lambda_1 +\dots x_n\otimes \lambda_n$  with $x_i\in E$ and $\lambda_i\in \K$, such that $n$ is minimal.
  After scalar multiplication we can assume that $\lambda_1\in \F$.
  
  Now $\hat{g}(r,y\otimes 1)=0$ for all $y\in E$,
  but, then also
  $\hat{g}(r^{\sigma},y\otimes 1)=0$ for all $y\in E$
  and hence $r-r^{\sigma}=x_2\otimes(\lambda_2-\lambda_2^{\sigma})+\dots + x_n\otimes(\lambda_n-\lambda_n^{\sigma})$
  is an element of the radical,
  which by assumption on minimality of $n$ is zero.
  But that would imply that $r$ is
  fixed by $\sigma$ and contradicts that the radical
  of $g$ is trivial.
\end{proof}

\begin{cor}\label{hatsg simple}
  If the characteristic of $\mathbb{F}$ is not $2$, then 
the Lie algebra $\hat{\g}$ is simple.
\end{cor}

\begin{proof}
  By Lemma \ref{radical trivial} the radical of $\hat{g}$ is trivial, and
  $\hat{g}$ is non-degenerate.
  
  Let $\mathfrak{i}$ be a nontrivial ideal of $\hat{\g}$. Then there is an extremal element
  $x\in\hat{E}$, the set of extremal elements of $\hat{\g}$, such that $\hat{g}(x,i)\neq 0$ for some $i\in \mathfrak{i}$.
  But then $[x,[x,i]]\in \mathfrak{i}$ and hence $ x\in \mathfrak{i}$.

  As, by Proposition \ref{linear span}, the set $E$ spans $\g$ and hence also $\hat{\g}$, there is an element $y\in E$ with $\hat{g}(x,y)\neq 0$, from which it follows that $[y,[y,x]]$ and thus also
  $y$ is in $\mathfrak{i}$. But then as $\Gamma_{\mathfrak{sl}_2}$ is connected by Theorem \ref{ci}, we also find $E$ to be contained in $\mathfrak{i}$, proving
  $\mathfrak{i}$ to be $\hat{\g}$ and $\hat{\g}$ to be simple.
\end{proof}

\section{Infinitesimal transvections}
\label{section:transvections}

Let $V$ be a vector space over a field $\mathbb{F}$ and let $v\in V$ and $\phi\in V^*$ be nonzero.
Then the linear map
$t_{v,\phi}:V\rightarrow V$
with $t_{v,\phi}(w)=\phi(w)v$ for all $w\in V$
is called an \emph{infinitesimal transvection} if $\phi(v)=0$ and
an \emph{infinitesimal reflection} if $\phi(v)\neq 0$.

Let $V$ be a vector space over a field $\F$
and let $\Pi\subseteq V^*$ be a subspace of the dual $V^*$ of $V$.
We assume that $\ann_V(\Pi)$ is $\{0\}$.
For finite dimensional $V$ we have $\Pi=V^*$ and it is well known that
the infinitesimal transvections and reflections generate the \emph{ general linear Lie algebra}
$\mathfrak{gl}(V)$, while the infinitesimal transvections generate $\mathfrak{sl}(V)$, the \emph{special linear Lie algebra}.
For possibly infinite dimensional $V$ and $\Pi=V^*$, these elements generate the so called \emph{finitary general linear Lie algebra} $\mathfrak{fgl}(V)$ and the \emph{finitary special linear Lie algebra} $\mathfrak{fsl}(V)$, respectively.
(A linear map of $V$ to itself is called \emph{finitary} if its kernel has finite codimension.)

For infinite dimensional $V$  the space $\Pi$ can be a proper
subspace of $V^*$ and we define (also in the finite dimensional case) the
Lie subalgebra of $\mathfrak{gl}(V)$ generated by the
infinitesimal transvections and reflections $t_{v,\phi}$ with $v\in V$ and $\phi\in \Pi$ to be $\mathfrak{fgl}(V,\Pi)$ and
by the
infinitesimal transvections $t_{v,\phi}$ with $v\in V$, $\phi\in \Pi$ and $\phi(v)=0$ to be  $\mathfrak{fsl}(V,\Pi)$.
Notice that this latter Lie algebra indeed consists of traceless finitary  linear transformations.

Let $t_{v,\phi}$ be an infinitesimal transvection. Then for each linear map $y:V\rightarrow V$
  we have for all $w\in V$:
  $$\begin{array}{ll}
    [t_{v,\phi},[t_{v,\phi},y]](w)
    &=(t_{v,\phi}(t_{v,\phi}y-yt_{v,\phi})-(t_{v,\phi}y-yt_{v,\phi})t_{v,\phi})(w)\\
    &=t_{v,\phi}(t_{v,\phi}(y(w)))-t_{v,\phi}(y(\phi(w)v))\\
    &\ -t_{v,\phi}(y(\phi(w)v))+yt_{v,\phi}(\phi(w)v)\\
    &=-2\phi(y(v))\phi(w)v\\
    &=-2\phi(y(v))t_{v,\phi}(w).\\
    \end{array}
  $$
Thus we find $$[t_{v,\phi},[t_{v,\phi},y]]=-2\phi(y(v))t_{v,\phi}.$$

 This implies that,
 at least for fields $\F$ of characteristic different from $2$,
 the infinitesimal transvections $t_{v,\phi}$, with $v\in V$ and $\phi\in \Pi$,
 are extremal in $\mathfrak{fsl}(V,\Pi)$.
 It is straightforward, but tedious, to check the Premet identities.
 So, also in even characteristic, the infinitesimal transvections are extremal.
 By the assumption that $\ann_V(\Pi)=\{0\}$ we find these elements not to be sandwiches.
 
 The Lie algebras $\mathfrak{pfgl}(V,\Pi)$ and $\mathfrak{pfsl}(V,\Pi)$
 are obtained by factoring out the center (which might be trivial) of $\mathfrak{fgl}(V,\Pi)$
 and $\mathfrak{fsl}(V,\Pi)$.
 
As a vector space $\mathfrak{fgl}(V,\Pi)$ is isomorphic with $V\otimes \Pi$.
An isomorphism is provided by the linear expansion of the map $t_{v,\phi}\mapsto v\otimes \phi$, where $v\in V$ and $\phi\in \Pi$.
The Lie product then translates to  
$$[v\otimes \phi,w\otimes \psi]=\phi(w)v\otimes \psi -\psi(v) w\otimes \phi$$
for pure tensors $v\otimes \phi$ and $w\otimes \psi\in V\otimes \Pi$.

We will  identify $\mathfrak{fgl}(V,\Pi)$ (as well as various of its subalgebras) with the (subalgebras of the) Lie algebra $\mathfrak{g}(V,\Pi)$ defined by this product on $V\otimes \Pi$.  
In particular, we identify $\mathfrak{fsl}(V,\Pi)$ with the subalgebra
$\mathfrak{s}(V,\Pi)$ generated by the pure tensors $v\otimes \phi\in V\otimes \Pi$,  with
$\phi(v)=0$.

Now suppose $h:V\times V\rightarrow\mathbb{F}$ is a non-degenerate skew-Hermitian form with respect to a field automorphism $\sigma$ of order $\leq 2$.
So,  $h$ is linear in the second coordinate,
$$h(v,w)=-h(w,v)^\sigma$$
for all $v,w\in V$,
and $h(v,w)=0$ for all $w\in W$ implies $v=0$.

If $v\in V$, then denote by $h(v,\cdot)$ the linear map $w\in V\mapsto h(v,w)$.
Then 

$$\begin{array}{ll}
  [v\otimes h(v,\cdot),w\otimes h(w,\cdot)]&=h(v,w)v\otimes h(w,\cdot)-h(w,v) w\otimes h(v,\cdot)\\
    &=h(v,w)v\otimes h(w,\cdot)+h(v,w)^\sigma w\otimes h(v,\cdot)\\
    &=(h(v,w)v+w)\otimes h(h(v,w)v+w,\cdot)\\
    &\phantom{=}-h(v,w)h(v,w)^\sigma v\otimes h(v,\cdot) -w\otimes h(w,\cdot). 
\end{array}$$

So, the elements $v\otimes h(v,\cdot)$ generate a Lie subalgebra of $\mathfrak{g}(V,V^*)$
over the field $\mathbb{F}_\sigma=\{\lambda\in \mathbb{F}\mid \lambda^\sigma=\lambda\}$.
This subalgebra is denoted by $\mathfrak{g}_h(V,V^*)$.
The subalgebra of $\mathfrak{s}(V,V^*)$ generated by the elements $v\otimes h(v,\cdot)$ with $h(v,v)=0$ is
denoted by  $\mathfrak{s}_h(V,V^*)$.

If $h$ is alternating (i.e., $\sigma$ is the identity), the algebra $\mathfrak{s}_h(V,V^*)$ is isomorphic to the finitary symplectic Lie algebra $\mathfrak{fsp}(V,h):=\langle x\in \mathfrak{sl}(V,V^*)\mid h(x(v),w)+h(v,x(w))=0$ for all $v,w\in V\}$, see \cite{CFsymp17},
while for proper skew-Hermitian forms (i.e., $\sigma$ has order $2$), the Lie algebra  $\mathfrak{s}_h(V,V^*)$ over $\mathbb{F}_\sigma:=\{\lambda\in \mathbb{F}\mid \lambda^\sigma=\lambda\}$ is isomorphic to a subalgebra of
the finitary special unitary Lie algebra
$\mathfrak{fsu}(V,h)=\{x\in \mathfrak{fsl}(V,V^*)\mid h(x(v),w)+h(v,x(w))=0$ for all $v,w\in V\}$.
If $V$ is generated by its isotropic vectors, i.e. $v\in V$ with $h(v,v)=0$ , we find this subalgebra to be the full finitary special unitary Lie algebra, as will be shown below.

Notice that these Lie algebras are obtained as centralizers of the automorphism of $\mathfrak{s}(V,\{h(v,\cdot)\mid v\in V\})$, defined as the semi-linear extension of the map

$$(\lambda (v\otimes h(w,\cdot))\mapsto \lambda^\sigma w\otimes h(v,\cdot)$$
where $v,w\in V$ with $h(v,w)=0$ and $\lambda\in \mathbb{F}$.

Versions of the  following results can also be found in \cite{fle15}.

\begin{prop}\label{generation}
Let $h$ be a non-degenerate skew-Hermitian form on the vector space $V$ over the field  $\mathbb{F}$
with respect to the field automorphism $\sigma$ of order $2$.
The Lie algebra $\mathfrak{fu}(V,h)$ over $\mathbb{F}_\sigma$
is linearly spanned by its  infinitesimal reflections and transvections.
\end{prop}

\begin{proof}
  
  We first consider the case where $\dim(V)=n< \infty$.
In this case it is well known that $\mathfrak{u}(V,h)$ has dimension $n^2$.

Let $v_1,\dots, v_n$ be a basis of $V$ such that the matrix of the form
$h$ with respect to this basis
is diagonal
$H=\mathrm{diag}(\lambda_1,\cdots,\lambda_n)$, 
where $\lambda_i\in \MF$ with
$\lambda_i^\sigma=-\lambda_i$.

Now consider the elements $${v_i\otimes h(v_i,\cdot)},\ {(v_j+v_l)\otimes h(v_j+v_l,\cdot)},\ \mathrm{and}\ {(v_j+\mu v_l)\otimes h(v_j+\mu v_l,\cdot)},$$
where $1\leq i,j,l \leq n$, $j<l$ and $\mu\in\MF$ a fixed element  with $\mu^\sigma\neq \mu$.

We can verify that these
$n^2$ elements form an independent set in $\mathfrak{u}(V,h)$.
So, as the dimension of  $\mathfrak{u}(V,h)$ equals $n^2$,
we have shown that $\mathfrak{u}(V,h)$ is spanned by its
infinitesimal transvections and reflections.

From the finite dimensional case, the finitary case follows easily.
\end{proof}

We now focus on subalgebras of the unitary Lie algebra generated by
infinitesimal transvections.
The existence of such transvections implies the existence of isotropic vectors
in $V$, i.e. $0\neq v\in V$  
with $h(v,v)= 0$.
Two  isotropic vectors $v,w\in V$ with $h(v,w)\neq 0$ span a so-called \emph{hyperbolic $2$-space}.
Isotropic vectors do exist, if the form $h$ is  \emph{trace valued}, that is
$h(v,v)\in \{\lambda-\lambda^\sigma\mid \lambda\in \mathbb{F}\}$ for all $v\in V$.
We have the following characterization of trace-valued forms:

\begin{prop}\label{trace valued} Suppose $V$ is an $\mathbb{F}$-vector space and $\sigma$ a field
  automorphism of $\mathbb{F}$ of order $2$.
  Let $h$ be a non-degenerate skew-Hermitian form with respect to the
  involution $\sigma$. Suppose $V$ contains a nonzero isotropic vector.

  Then the form $h$ is trace-valued if and only if $V$  is generated by its isotropic vectors,
if and only if for each isotropic $0\neq v\in V$
    and arbitrary $w\in V$ with $h(v,w)\neq 0$, the subspace $\langle v,w\rangle$ is hyperbolic.
\end{prop}

\begin{proof}
See \cite[Theorem 10.1.3]{BC13}.
\end{proof}

\begin{prop}\label{transvection generation}
  Suppose $V$ is an $\mathbb{F}$-vector space of dimension at least $2$
  and $\sigma$ a field
  automorphism of $\mathbb{F}$ of order $2$.
  Let $h$ be a non-degenerate skew-Hermitian form with respect to the
  involution $\sigma$ such that $V$ contains a nonzero isotropic vector.

  Then $\mathfrak{fsu}(V,h)$ is linearly spanned by its infinitesimal transvections if and only if $h$ is trace-valued.
\end{prop}

\begin{proof}
  Suppose $h$ is trace-valued.
  Then, by Proposition \ref{trace valued}, we can assume that $V$ is spanned by its isotropic vectors. 
We will show that  $\mathfrak{fsu}(V,h)$ is linearly spanned by its
infinitesimal transvections.

In view of the above result, Proposition \ref{generation}, and the fact that  $\mathfrak{fsu}(V,h)$
has codimension $1$ in $\mathfrak{fu}(V,h)$, it suffices to prove that
$\mathfrak{fu}(V,h)$
can be spanned by all its infinitesimal transvections and a unique reflection.

This is clearly true in the case where $(V,h)$ is a hyperbolic $2$-space
over the field $\mathbb{F}$.
Indeed, if $v_1,v_2$ is a hyperbolic basis of $V$, i.e. $h(v_1,v_1)=h(v_2,v_2)=0$
and $h(v_1,v_2)\neq 0$,
then $${v_1\otimes h(v_1,\cdot)},\ {v_2\otimes h(v_2,\cdot)},\ \mathrm{and}\  {v_1+v_2\otimes h(v_1+v_2,\cdot)}$$ span $\mathfrak{su}(V,h)$
and together with any infinitesimal reflection they span
$\mathfrak{u}(V,h)$.

So, to prove in general that
$\mathfrak{fu}(V,h)$ can be spanned by its infinitesimal transvections
together with one reflection, it suffices to prove connectedness of
the graph $\Delta$ on the anisotropic points of $V$, where two such points
are adjacent if and only if they span a hyperbolic line (i.e., $2$-space).
For then, every infinitesimal reflection is contained in the span of the
infinitesimal transvections and the unique reflection.

Suppose that $\dim(V)$ is at least $2$. 
As $(V,h)$ is non-degenerate and  spanned by its  isotropic vectors, there are isotropic vectors
$v,w\neq 0$  with $h(v,w)\neq 0$. After replacing $w$ with $h(v,w)^{-1}w$,
we can assume $h(v,w)=1$.
The $2$-space $\langle v,w\rangle$ is hyperbolic.

Let $u$ be an anisotropic vector  in $V$.
We will prove that $\langle u\rangle $ is in the same connected component
of $\Delta$ as some anisotropic point on $\langle v,w\rangle$.

First assume that $h(u,v)=h(u,w)=0$.
Then, for  $\mu\in \mathbb{F}$ satisfying $\mu-\mu^\sigma=-h(u,u)$
we have $h(u+v+\mu w,u+v+\mu w)=h(u,u)-\mu^\sigma +\mu=0$
and $h(v+\mu w,v+\mu w)=\mu-\mu^\sigma=-h(u,u)\neq 0$.
Notice, such $\mu$ exists, as $h$ is trace-valued.
Moreover, as $h(u,u+v+\mu w)=h(u,u)\neq 0$, we
find the $2$-space $\langle u, v+\mu w\rangle$ to be  hyperbolic, see \ref{trace valued}.
This implies that $\langle u\rangle$ is adjacent to $\langle v+\mu w\rangle$
in $\Delta$.

Now assume that $h(u,v)\neq 0$ but $h(u,w)=0$.
After scaling $v$ (and $w$), we can assume that $h(u,u)-h(u,v)\neq 0$.
Let $\mu\in \mathbb{F}$ with
$ h(u,u)+(\mu-h(u,v))^\sigma-(\mu-h(u,v))=0$.
Then $h(u+v+\mu w,u+v+\mu w)=
h(u,u)+h(u,v)-h(v,u)+\mu^\sigma-\mu=
h(u,u)+(\mu-h(u,v))^\sigma-(\mu-h(u,v))=0$.
Again, such element $\mu$ exists as $h$ is trace valued.
Moreover, $h(u,u+v+\mu w)=h(u,u)-h(u,v)\neq 0$.
So, $\langle u, u+v+\mu w\rangle$ is hyperbolic and meets $\langle v,w\rangle$
in $v+\mu w$, which is anisotropic, as $h(v+\mu w,v+\mu w)=\mu-\mu^\sigma=h(u,u)$.
We find that $\langle u\rangle$ is adjacent in $\Delta$ to an anisotropic point on $\langle v,w\rangle$.

Finally assume that $h(u,v)\neq 0$ and $h(u,w)\neq 0$.
Let $u'\in \langle u,v\rangle$ be perpendicular to $w$.
If $u'$ is anisotropic, then
$\langle u\rangle$ is adjacent to
$\langle u'\rangle$, and the latter is, by the above, adjacent to some anisotropic point in $\langle v,w\rangle$.
Thus, assume  $u'$ is isotropic. Then all nonzero vectors in
$\langle u',w\rangle$ are isotropic.
As $h(u,u')\neq 0$ we find that $u$ is perpendicular to  one
point on $\langle u',w\rangle$. So, all $2$-spaces on $\langle u\rangle$ inside $\langle u,v,w\rangle$,  except for  one, are hyperbolic.
This clearly implies that there is at least one hyperbolic line on
$\langle u\rangle$ meeting  $\langle v,w\rangle$ in a anisotropic point,
and that $\langle u\rangle$ is adjacent in $\Delta$ to an anisotropic point
on $\langle v,w\rangle$, finishing the proof that $\mathfrak{fsu}(V,h)$ is spanned by its transvections.

Now assume that the form $h$ is not trace-valued and thus, by Proposition \ref{trace valued}, that the subspace $V_0$ of $V$ spanned by the isotropic vectors
is not $V$.
Then we can find a vector $u\in V\setminus V_0$.
Clearly $u\otimes h(u,\cdot)$ is  not
in the subspace of $\mathfrak{fu}(V,h)$
spanned by the infinitesimal transvections and all elements
$v\otimes h(v,\cdot)$, with $v\in V_0$, including  infinitesimal reflections.
So,  the subalgebra spanned by the infinitesimal transvections has codimension at least $2$ in $\mathfrak{fu}(V,h)$ and hence is properly contained in $\mathfrak{fsu}(V,h)$, which has codimension $1$.
This finishes the proof of the proposition.
\end{proof}

Let $h$ be a skew-Hermitian
form on the vector space $V$
  defined over the field  $\mathbb{F}$
  with respect to the field automorphism $\sigma$ of order $2$.
  Assume that $h$ is trace-valued and $V$ contains an isotropic vector.

  Let $\g$ be the Lie algebra $\mathfrak{fsu}(V,h)$, which, by Proposition \ref{transvection generation}, is spanned by its infinitesimal transvections.
  Then  $\hat{\g}:=\g\otimes_{\mathbb{F}_\sigma}\mathbb{F}$ is
  isomorphic to a subalgebra of $\mathfrak{fsl}(V,\Pi)$ where $\Pi$ is the subspace of $V^*$
  spanned by the elements $h(v,\cdot)$ with $v\in V$ isotropic.
  If $V$ is finite dimensional, then, as $\hat{\g}$ and $\mathfrak{fsl}(V,\Pi)$
  have the same dimension, they are even isomorphic.
  This also holds true for infinite dimensional $V$, as we can view
  both Lie algebras as limits of isomorphic finite dimensional ones.
  See also the next section.

  It is well known that the Lie algebra $\mathfrak{pfsl}(V,\Pi)$ is  simple.
  We use this to find the structure of $\g$.

  As both $\g$ and $\hat{\g}$ are generated by their pure extremal elements, they are both equipped with an
  extremal form, say $g$ and $\hat{g}$, respectively.

  \begin{lemma}
The radical of $g$ equals the center of $\g$, which is of dimension at most $1$.
  \end{lemma}

  \begin{proof}
    Let $R=\rad(g)$. Then $\hat{R}:=R\otimes \mathbb{F}$ is contained in the radical
    $\rad(\hat{g})$, which is a proper ideal of $\hat{\g}$.
    As $\hat{\g}$ modulo its center is simple, we find  $\hat{R}$ to be contained in the center of $\hat{\g}$. But then
    $R$ is contained in the center of $\g$.

    Since the center of $\hat{\g}$ is at most one-dimensional, so is the
    center of $\g$.
    \end{proof}

\begin{prop}
  $\g/\rad(g)$ is simple.
\end{prop}

\begin{proof}
  Let $\mathfrak{i}$ be a nontrivial ideal of $\g$.
  
  If  $\mathfrak{i}\leq \rad(g)$, then it is central.
  If $\mathfrak{i}$ contains an element $i$ not in $\rad(g)$, we find an
  infinitesimal transvection $t\in \g$ with $0\neq t_0=[t,[t,i]]\in \mathfrak{i}$.
  So $t$ is in $\mathfrak{i}$.
  But then for each infinitesimal transvection $t'$ with $[t',t_0]\neq 0$ we find $t'\in \mathfrak{i}$.
  Repeating this argument we find all infinitesimal transvections
  in $\mathfrak{i}$ and $\mathfrak{i}=\g$.
\end{proof}

\section{A geometric characterization of the special linear Lie algebra}
\label{section:sl}

  In this section we  prove Theorem \ref{slthm}
 and show that  a simple Lie algebra  over a field $\F$ generated by its set of pure extremal elements
  and its extremal geometry isomorphic to $\Gamma(V,\Pi)$,
  where $V$ is an $\F$-vector space of dimension at least $3$ and $\Pi$ a subspace of $V^*$
  with $\ann_V(\Pi)=\{0\}$, is isomorphic to $\mathfrak{(p)fsl}(V,\Pi)$.

  If $V$ is finite dimensional, then the theorem follows from \cite{CRS} ( see also \cite{Rob12}).
  Indeed, using \cite{CRS} we obtain the following theorem.
  (Here we say that a subset $\E_0$ of $\E$ is \emph{closed under the action} of a subgroup $G$ of $\mathrm{Aut}(\g)$ if and only if $x^g\in \E_0$ for all $x\in \E_0$ and $g\in G$.)

  \begin{thm}[Cuypers, Roberts, and Shpectorov]\label{finite sl thm}
    Let $V$ be a vector space over a field $\mathbb{F}$ of finite  dimension at least $3$ and $\g$ a Lie algebra generated by a subset
    $\E_0$  of the set of  pure extremal points $\E$ such that $\E_0$ is a subspace
    of the extremal geometry $\Gamma(\g)$,
    isomorphic to $\Gamma(V,V^*)$, and closed under the action of $\mathrm{Exp}(x)$ for each $x\in \E_0$.

    Then $\g$ is isomorphic to $\mathfrak{sl}(V)$ or its central quotient $\mathfrak{psl}(V)$.
  \end{thm}

  \begin{proof}
    Although in the statements of the main results of \cite{CRS} it is required that $\E_0=\E$  and that the Lie algebra is simple,
    this is actually not needed in the proof.

    Indeed, we can take a basis $v_1,\dots, v_n$ of $V$ and consider the extremal elements $x_{ij}$ with  $\langle x_{ij}\rangle\in \E_0$, where $i\neq j$,
    corresponding to the point-hyperplane pairs $(\langle v_i\rangle, \langle v_k\mid k\neq j\rangle)$.
    As is shown in Section 5 of \cite{CRS}, we can scale these elements in such a way that they together with the elements $h_{ij}=-h_{ji}=[x_{ij},x_{ji}]$ form a Chevalley spanning set as described in Section 4 of \cite{CRS}.
      But then, using the group $G=\langle \mathrm{Exp}(x)\mid x\in \E_0\rangle$ it readily follows by arguments as in  Section 6 and 7 of \cite{CRS} that the Lie algebra generated by this Chevalley spanning set
    contains $\E_0$ and hence equals $\g$. But then $\g$ is a Chevalley Lie algebra of type $\mathrm{A}_{n-1}$ and hence isomorphic to
    $\mathfrak{(p)sl}(V)$.
    \end{proof}

So, for the remainder of this section we can and will 
assume that $V$ is infinite dimensional.

Our proof of Theorem \ref{slthm}  will rely on Theorem \ref{finite sl thm}.
We will approach an unknown infinite dimensional Lie algebra satisfying the hypothesis of the
theorem by a collection of finite dimensional
  subalgebras isomorphic to $\mathfrak{sl}(U)$ for vector spaces $U$ of finite dimension.

  We start with the necessary definitions.

  \bigskip

  A {\em directed set} $(\I,\sqsubseteq)$ is a set  $\I$ with a partial order $\sqsubseteq$ such that
  \begin{itemize}

  \item
    for all $i,j\in \I$ there exists an element $k$ in  $\I$ such that
    $i,j\sqsubseteq k$;
  \item for every $k\in \I$ any chain
    $i_1\sqsubseteq i_2\sqsubseteq \dots$
    of distinct elements $i_j\sqsubseteq k$ from $\I$  has finite length.
  \end{itemize}

Suppose $\g$ is a Lie algebra and 
$(\I,\sqsubseteq)$  a directed set.
Then a collection  $(\g_i)_{i\in \I}$  of Lie subalgebras of
$\g$ is called a {\em local system} of $\g$ with respect to $(\I,\sqsubseteq)$ if the following holds:

\begin{itemize}
\item if $i,j\in \I$ with $i\sqsubseteq j$, then $\g_i\subseteq \g_j$;
  \item $\g=\displaystyle\bigcup_{i\in \I}\g_i$
\end{itemize}

If $\g$ contains a local system $(\g_i)_{i\in \I}$, then $\g$ is uniquely determined
by this local system as it is isomorphic to the direct limit
$$\displaystyle \lim_{i\in \I}\g_i.$$

This implies the following well known result.

\begin{prop}\label{directed set iso}
  Let $(\I,\sqsubseteq)$ be a directed set.
  Suppose both $\g$ and $\h$ are  Lie algebras with local systems
  $(\g_i)_{i\in \I}$ and $(\h_i)_{i\in \I}$ for $(\I,\sqsubseteq)$, respectively, such that
  for each $i\in \I$ there exists an isomorphism
  $$\phi_i:\g_i\rightarrow \h_i$$
  with the property that for $j,k\in \I$ with $j\sqsubseteq k$ we have
  $$\phi_j=\phi_k|_{\g_j}.$$
  Then $\g$ and $\h$ are isomorphic.
\end{prop}

We now start with the proof of Theorem \ref{slthm}.
So, suppose $V$ is an infinite dimensional vector space over the field $\F$
and $\Pi$ a subspace of $V^*$ with $\ann_V(\Pi)=\{0\}$.

Then let $\I$ be the set of all pairs $(U,\Phi)$ where
\begin{itemize}
  \item $U$ is a subspace of $V$;
  \item $\Phi$ is a subspace of $\Pi$;
  \item $U$ and $\Phi$ have the same finite dimension, which is at least $3$ and \emph{not} divisible by the characteristic of the field $\F$;
  \item $\ann_U(\Phi)=\{0\}$.
\end{itemize}

For $(U,\Phi), (U',\Phi')\in \I$ we define $\sqsubseteq$ by the following:
$$(U,\Phi)\sqsubseteq (U',\Phi')\Leftrightarrow U\subseteq U' \ \mathrm{and}
 \  \Phi\subseteq \Phi'.$$

\begin{lemma}
$(\I,\sqsubseteq)$ is a directed set.
\end{lemma}

\begin{proof}
  Let $(U_1,\Phi_1)$ and $(U_2,\Phi_2)$ be elements in $\mathcal{I}$.
  Now let $U=U_1+U_2$ and $\Phi=\Phi_1+\Phi_2$.
  
  Suppose $\phi_1,\dots,\phi_m$ form a basis for $\Phi$.
  Then pick vectors $u_1,\dots, u_m\in V$ with $\phi_i(u_j)=\delta_{ij}$
  for all $1\leq i,j \leq n$ and consider $\hat{U}=\langle u_1,\dots,u_m\rangle +U$.
  Let $u_{m+1},\dots, u_n$ be vectors such that $u_1,\dots,u_n$
  form a basis for $\hat{U}$.
  Then we can find $\phi_{m+1},\dots, \phi_n\in \Pi$ such that $\phi_i(u_j)=\delta_{ij}$  for all $1\leq i,j\leq n$.
  Now define $\hat{\Pi}$ to be $\langle \phi_1,\dots,\phi_n\rangle$.
  Then $\dim(\hat{U})=\dim(\hat{\Phi})=n$ and $\ann_{\hat{U}}(\hat{\Phi})=\{0\}$.

  If the dimension $n$ is divisible by the characteristic of $\F$,
  then extend $\hat{U}$ and $\hat{\Phi}$ with  elements $u$ and $\phi$,
  respectively, satisfying $\phi(u_i)=0=\phi_i(u)$ for all $1\leq i\leq n$ and
  $\phi(u)=1$.
  In any case we find an element $(\hat{U},\hat{\Phi})$ of $\mathcal{I}$
  with $(U_1,\Phi_1),(U_2,\Phi_2)\sqsubseteq (\hat{U},\hat{\Phi})$.

  Since every pair $(U,\Phi)\in \I$ consists of finite dimensional spaces, we find $(\I,\sqsubseteq)$ to be a directed set.
\end{proof}

We now construct a local system for $\mathfrak{fsl}(V,\Pi)$.

For $I=(U,\Phi)\in \I$  define
$$\mathfrak{sl}(I):=\langle
t_{v,\phi}\mid 0\neq v\in U, 0\neq \phi\in \Phi \ \mathrm{and} \ \phi(v)=0\rangle.$$
Then $\mathfrak{sl}(I)$ is a Lie subalgebra
of $\mathfrak{fsl}(V,\Pi)$ which is isomorphic to $\mathfrak{sl}(U)$, and,
as a consequence of the  extra condition that the dimension of $U$ and $\Phi$
is not divisible by the characteristic of $\F$,
has a trivial center.

\begin{prop}
   The collection of Lie subalgebras $\mathfrak{sl}(I)$ for $I\in \I$ forms
   a local system for $\mathfrak{fsl}(V,\Pi)$.
\end{prop}

Now consider a simple Lie algebra $\g$ generated by its set $E$ of pure extremal elements and assume its extremal geometry $\Gamma(\g)$ is isomorphic
to $\Gamma(V,\Pi)$.
We also construct a local system for $\g$.
We identify $\Gamma(\g)$ with $\Gamma=\Gamma(V,\Pi)$ and for each $I=(U,\Phi)\in \mathcal{I}$ define $\Gamma_I$ to be the subspace of $\Gamma$
consisting of all incident point-hyperplane pairs $(p,H)$ with $p\in \mathbb{P}(U)$ and  $H=\ker(\phi)$ for some $0\neq \phi\in \Phi$.

Moreover, we define $\g_{I}$ to be the Lie subalgebra generated by all extremal points in $\Gamma_I$. So,
$$\g_I:=\langle x\mid x\in \Gamma_I\rangle.$$

\begin{lemma}
  Let $I\in \I$. Then $\E\cap \g_I=\Gamma_I$.
\end{lemma}

\begin{proof}
  Let $I=(U,\Phi)$ be an element from $\I$ and take $U'=\ann_V(\Phi)=\{u\in V\mid \phi(u)=0$ for all $\phi\in \Phi\}$
  and $\Phi'=\ann_\Pi(U):=\{\phi\in\Pi \mid \phi(u)=0 $ for all $u\in U\}$.
  Then all extremal points $x$ from $\g$ which correspond to points $(p,H)$ from $\Gamma$ with $p$ in $\mathbb{P}(U')$ and 
  $H=\ker(\phi')$ for some $\phi'\in \Phi'$ centralize
  all elements from $\Gamma_I$ and hence $\g_I$.

  Suppose  an extremal point $y$ of $\g_I$
  corresponds to the point-hyperplane pair $(q,K)$ of $\Gamma$. Then $y$  is centralized by all such
  extremal points $x$, and we find $q\leq \mathrm{Ann}_V(\Phi')=U$, and
  $K=\ker(\phi)$ for some $\phi\in \ann_\Pi(U')=\Phi$.
  In particular, $y\in \Gamma_I$.
  \end{proof}

\begin{prop}
   The collection of Lie subalgebras $\g_I$ for $I\in \I$ forms
   a local system for $\g$.
\end{prop}

\begin{prop}
  The two local systems $(\mathfrak{sl}(I))_{I\in \mathcal{I}}$ and
  $(\g_I)_{I\in \mathcal{I}}$ are isomorphic.
\end{prop}

\begin{proof}    
    For each $I=(U,\Phi)\in \I$ we find, as we can apply Theorem \ref{finite sl thm}, that the Lie algebra $\g_I$ is classical Lie algebra of type $A_{n-1}$, where $n=\dim(U)$.
    In particular, $\g_I$ is isomorphic to  $\mathfrak{sl}(U)$ or its central quotient $\mathfrak{psl}(U)$.
    As $\dim(U)$ is not divisible by the characteristic of $\F$, we find
    $\g_I$ to be isomorphic to $\mathfrak{sl}(U)$.
    Moreover, the isomorphism can be chosen in such a way that
    it induces the identity on $\Gamma_I$.

    It remains to show that for $I\sqsubseteq J$
    the restriction of $\phi_J$ to $\g_I$ equals $\phi_I$.

    So, suppose $I=(U,\Phi)\sqsubseteq J$.
   Then the map $\psi=\phi_I(\phi_J|_{\g_I})^{-1}$ is an automorphism of $\mathfrak{sl}(I)\simeq \mathfrak{sl}(U)$ fixing all extremal points of $\mathfrak{sl}(I)$.
    Let $x,y$ be two collinear points of $\Gamma_I$,
    then $\psi$ fixes all points on this line and hence is a scalar multiplication on this $2$-dimensional subspace of $\mathfrak{sl}(I)$.
    By connectedness of $\Gamma_I$ we find that $\psi$ is a scalar multiplication on the vector space
    $\mathfrak{sl}(I)$, say with scalar $\lambda\in \F^*$.
    But as $\psi$ also respects the Lie product $[\cdot,\cdot]$ of $\mathfrak{sl}(I)$, we find for any two  elements $a,b\in \mathfrak{sl}(I)$
    that $[\psi(a),\psi(b)]=[\lambda a, \lambda b]=\lambda^2[a,b]=\psi([a,b])=\lambda [a,b]$.
    In particular, by choosing non commuting $a,b\in \mathfrak{sl}(I)$ we find $\lambda^2=\lambda$, and hence $\lambda =1$.

    So indeed, $\phi_J|_{\g_I}=\phi_I$ and we have found the two local systems to be isomorphic.
\end{proof}

By Proposition \ref{directed set iso} we immediately have the following result, which together with the finite dimensional case, settles Theorem \ref{slthm}.

\begin{thm}\label{sl thm}
  The Lie algebra $\g$ is isomorphic to $\mathfrak{fsl}(V,\Pi)$.
\end{thm}

\section{A geometric characterization of the unitary  Lie algebra}
\label{section:u}
In this final section we focus on Theorem \ref{uthm}.
So, let $\g$ be a simple Lie algebra over the field $\F$ generated by its set $E$ of pure extremal elements, and suppose that the extremal geometry $\Gamma(\g)$ does not contain lines.

In this section we assume that there is a quadratic Galois extension $\K$ of $\F$ such that the set of pure extremal elements $\hat{E}$ of 
$\hat{\g}=\g\otimes_{\F}\K$ has an extremal geometry $\hat{\Gamma}:=\Gamma(\hat{\g})$
isomorphic to $\Gamma(V,\Pi)$ where $V$ is a $\K$-vector space and $\Pi$ a subspace of $V^*$ with $\ann_V(\Pi)=\{0\}$.
The Lie algebra $\hat{\g}$ contains pure extremal elements and,
by Theorem \ref{finite sl thm} and \ref{sl thm}, the Lie algebra $\hat{\g}$ is isomorphic to $\mathfrak{pfsl}(V,\Pi)$.
(Indeed, if the center is nontrivial, then it is fixed by $\sigma$ and we find  $\g$ to have a nontrivial center.) 

As already observed in Section \ref{section:extremal}, the nontrivial field automorphisms $\sigma$ of the extension induces an  automorphism, also denoted by 
$\sigma$ of $\hat{\g}$ acting on an element $x\otimes \lambda$, with $x\in \g$ and $\lambda\in \mathbb{K}$  as

$$(x\otimes \lambda)^{\sigma}=x\otimes \lambda^\sigma.$$

This automorphism also induces an automorphism, which we also denote by $\sigma$, on the extremal geometry $\hat{\Gamma}$.

Let $\hat{\E}$ be the set of extremal points of $\hat{\g}$.
Then, as we assume that $\hat{\Gamma}$ is isomorphic to
$\Gamma(V,\Pi)$ for some vector space $V$ and subspace $\Pi$ of $V^*$
with annihilator $\ann_V(\Pi)=\{0\}$,
we can identify the elements of $\hat{\E}$
with incident pairs $(p,H)$, where $p\in \mathbb{P}:=\mathbb{P}(V)$
and $H\in \mathbb{H}$, where $\mathbb{H}$ is the set of hyperplanes
$\ker(\phi)$, with $0\neq \phi\in \Pi$.

  Let $p\in \mathbb{P}$ and $H\in \mathbb{H}$.
  Then   $\overline{p}:=\{(p,K)\mid K\in \mathbb{H}$ with $p\in K\}$ and 
  $\overline{H}:=\{(q,H)\mid q\in \mathbb{P}$ with $q\in H\}$
  are maximal cliques of the collinearity graph of $\Gamma(V,\Pi)$.
  Let $\overline{\mathbb{P}}$ denote the set of cliques $\overline{p}$
  with $p\in \mathbb{P}$ and
  $\overline{\mathbb{H}}$ the set of cliques $\overline{H}$ with $H\in \mathbb{H}$.
    Finally set $\mathcal{C}$ to be the union of $\overline{\mathbb{P}}$ and
  $\overline{\mathbb{H}}$.

    \begin{prop}
      \begin{enumerate}[\rm (a)]
      \item
        
      Let $C$ be a clique in the collinearity graph of $\Gamma(V,\Pi)$ of size at least $2$, then
      $C$ is contained in a unique clique from $\mathcal{C}$.

    \item A point $p\in \mathbb{P}$ is incident with a hyperplane $H\in \mathbb{H}$ if and only if $\overline{p}\cap \overline{H}\neq \emptyset$.

    \item The elements from $\overline{\mathbb{P}}$ partition the point set of $\Gamma(V,\Pi)$.

    \item The elements from $\overline{\mathbb{H}}$ partition the point set of $\Gamma(V,\Pi)$.

      \end{enumerate}

  \end{prop}

\begin{proof}
Straightforward.
  \end{proof}

The set $\mathcal{C}$ is partitioned into the two parts $\overline{\mathbb{P}}\cup\overline{\mathbb{H}}$.
To be in the same part of the partition can be characterized as follows:

\begin{lemma}
  Two maximal cliques $C_1$ and $C_2$ are in the same part
  of the partition $\overline{\mathbb{P}}\cup\overline{\mathbb{H}}$ if and only if there is a third maximal clique $C_3$
with $$C_1\cap C_3\neq \emptyset\neq C_2\cap C_3.$$
\end{lemma}

\begin{proof}
Straightforward.
  \end{proof}

\begin{lemma}
  If $p\in \mathbb{P}$, then $\overline{p}^\sigma\in \overline{\mathbb{H}}$
  and if $H\in \mathbb{H}$, then $\overline{H}^\sigma\in \overline{\mathbb{P}}$.
  \end{lemma}

\begin{proof}
  Consider an element $x\in \hat{\Gamma}$ fixed by $\sigma$.
  Then $x=(p,H)$ for some point $p\in \mathbb{P}$ and hyperplane $H\in \mathbb{H}$.
  As $\sigma$ does not fix any line of $\hat{\Gamma}$, it must map
  $\overline{p}$ to $\overline{H}$.

  For every $q\in \mathbb{P}$ different from
  $p$, we find
  that
  $\overline{q}$ does not meet $\overline{p}$, but there is a maximal clique
  meeting both $\overline{p}$ and $\overline{q}$. 
  Since $\sigma$ is an automorphism of $\Gamma$, we find that
   $\overline{q}^\sigma$ does not meet $\overline{p}^\sigma$, but there is a maximal clique
  meeting both $\overline{p}^\sigma$ and $\overline{q}^\sigma$. So,
  by the above lemma, we also have $\overline{q}^\sigma\in \overline{\mathbb{H}}$.
The argument for hyperplanes is similar.
  
\end{proof}

\begin{lemma}
  If $p,q\in \mathbb{P}$, then $\overline{p}\cap \overline{q}^\sigma\neq \emptyset$
  if and only if $\overline{p}^\sigma\cap \overline{q}\neq \emptyset$
  \end{lemma}

\begin{proof}
  This is immediate since $\sigma$ is an involution.
\end{proof}

The above two lemmas imply that $\sigma$ induces a quasi-polarity on $\mathbb{P}$
mapping a point $p$ of $\mathbb{P}$ to the hyperplane $H$ with $\overline{H}=\overline{p}^\sigma$.
Again we use $\sigma$ to denote this quasi-polarity.
This  quasi-polarity $\sigma$ is non-degenerate, as each point is mapped
to a proper hyperplane in $\mathbb{H}$.

The pairs $(p,H)$ in $\E$ are precisely the pairs $(p,H)$ for which
$p$ is contained in $H= p^\sigma$.

The  quasi-polarity $\sigma$ can be realized by a reflexive sesquilinear form
$$h:V\times V\rightarrow \mathbb{K}.$$
That means, there are $\epsilon\in \mathbb{K}$ and $\tau$ in the automorphism group of $\mathbb{K}$ of order at most $2$
with $\epsilon^\tau=\epsilon ^{-1}$
such that $h$ is linear in the second coordinate and satisfies
$$h(v,w)=\epsilon h(w,v)^\tau$$
for all $v,w\in V$.

Due to Theorem \ref{slthm} and the above
we have the following  identifications:
$$x\in \hat{E} \leftrightarrow t_{v,\phi}=t_{v,h(w,\cdot)} \leftrightarrow v\otimes h(w,\cdot),$$
for some $w\in V$.

The polarity $\sigma$ determines the form $h$ up to proportionality (i.e. we can scale $h$ by multiplying it with a scalar).
But that implies that we can fix an isotropic vector $v_0\in V$ and scale $h$ in such way that
$$v_0\otimes h(v_0,\cdot)\in E.$$
Having done this we obtain the following.
\begin{lemma}
  We have $\epsilon=-1$  and $\tau=\sigma$,
  and $E$ consists of the elements $\alpha v\otimes h(v,\cdot)$
  with $0\neq v\in V$, $h(v,v)=0$ and $\alpha\in \mathbb{F}$.
\end{lemma}

\begin{proof}
  Let $x=v_0\otimes h(v_0,\cdot)$. Then we find $\mathrm{exp}(x,\lambda)$ with $\lambda\in \mathbb{F}$ to be an $\mathbb{F}$-linear automorphism of $\g$.
  In particular, if   $\langle w\otimes h(w,\cdot)\rangle$ is a $1$-space of $\hat{\g}$ containing an extremal element of $E$, we find
  that $\langle \mathrm{exp}(x,\lambda)(w\otimes h(w,\cdot))\rangle$ also contains elements of $E$.
  Take such an element with $h(v_0,w)\neq 0$.
  
  Then $$\mathrm{exp}(x,\lambda)(w\otimes h(w,\cdot))=(w+\lambda h(v_0,w) v_0)\otimes(h(w,\cdot)-\lambda h(w,v_0) h(v_0,\cdot)).$$
  As the latter spans a $1$-space containing an element of $E$, we do have that
  $h(w+\lambda h(v_0,w)v_0,\cdot)$ is a scalar multiple of $h(w,\cdot)-\lambda h(w,v_0)h(v_0,\cdot)$.

  But this implies that $$w-\lambda^\tau h(w,v_0)^\tau v_0=w-\epsilon\lambda^\tau h(v_0,w)^{\tau^2}v_0=w+\lambda h(v_0,w) v_0.$$
  In particular, we find $-\epsilon\lambda^\tau =\lambda$.
  Taking $\lambda$ to be $1$ yields $\epsilon=-1$ and then $\lambda^\tau=\lambda$ for all $\lambda\in \mathbb{F}$.
  So, $\tau=\sigma$ or the identity.

  If $\tau$ is the identity, then $h$ is a symplectic form.
  In particular, for each vector $0\neq v\in V$ we find
  $\langle v\rangle \subseteq \langle v\rangle^\sigma$.
  But then $\hat{\g}$ is contained in the symplectic subalgebra $\mathfrak{s}_h(V,V^*)$, and does not have
  an extremal geometry isomorphic to $\Gamma(V,\Pi)$.

So, $\tau=\sigma$.
  Now if $y=\alpha w\otimes h(w,\cdot)$ is in $E$, then
  consider
  $\exp(y,\lambda)x=(v_0+\lambda h(w,v_0)\alpha w)\otimes(h(v_0+\alpha^\sigma\lambda h(w,v_0)w,\cdot))\in E$.

  As above, this can only be in $E$ provided $\alpha^\sigma=\alpha$.
  But that proves the lemma.
\end{proof}

This results in the following theorem.

\begin{thm}
  $\g$ is isomorphic to $\mathfrak{s}_h(V,V^*)$ modulo its center, where
  $h$ is a non-degenerate skew-Hermitian form with $h(v,w)=-h(w,v)^\sigma$ for all $v,w\in V$.
\end{thm}

\begin{proof}
The lemma clearly implies that $\mathfrak{g}$ is, up to a center, isomorphic
to $\mathfrak{s}_h(V,V^*)$, for the skew-Hermitian form $h$.
\end{proof}

Since $\hat{\g}$ is isomorphic to $\mathfrak{s}(V,\Pi)$ we find that
the $V$ is generated by the  vectors $v\in V$ with $h(v,v)=0$.
This implies, by Proposition \ref{trace valued}, that the form $h$ is trace-valued. 
Moreover, by Proposition \ref{transvection generation}, the Lie algebra $\mathfrak{s}_h(V,V^*)$ modulo its center is isomorphic to $\mathfrak{pfsu}(V,h)$.
So, we have finished the proof of Theorem \ref{uthm}.

\bibliographystyle{plain}

\bibliography{marc.bib}

\end{document}